\documentclass[12pt,a4paper]{amsart}
\usepackage[utf8]{inputenc}	
\usepackage[T1]{fontenc}
\usepackage[in,headings]{fullpage}
\usepackage{amsbsy, amscd, amsfonts, amsmath, amsrefs, amssymb, amsthm}
\usepackage{graphicx}
\usepackage{float}
\usepackage{color}   
\usepackage[colorlinks=true,linkcolor=blue,citecolor=blue,urlcolor=blue]{hyperref}
\usepackage{comment}
\excludecomment{A}

\newtheorem{theorem}{Theorem}

\newtheorem{lemma}[theorem]{Lemma}

\theoremstyle{definition}

\newtheorem{remark}[theorem]{Remark}

\theoremstyle{remark}

\renewcommand{\Re}{\mathop{\mathrm{Re}}\nolimits}

\newfont{\cmbsy}{cmbsy10}
\newfont{\cmmib}{cmmib10}

\begin{document}

\title{On the approximation of the zeta function by Dirichlet polynomials}
\author[Arias de Reyna]{J. Arias de Reyna}
\address{%
Universidad de Sevilla \\ 
Facultad de Matem\'aticas \\ 
c/Tarfia, sn \\ 
41012-Sevilla \\ 
Spain.} 




\email{arias@us.es, ariasdereyna1947@gmail.com}


\begin{abstract}
We prove that for $s=\sigma+it$ with $\sigma\ge0$ and $0<t\le x$, we have \[\zeta(s)=\sum_{n\le x}n^{-s}+\frac{x^{1-s}}{(s-1)}+\Theta\frac{29}{14} x^{-\sigma},\qquad \frac{29}{14}=2.07142\dots\] where $\Theta$ is a complex number with $|\Theta|\le1$. This improves Theorem 4.11 of Titchmarsh. 
\end{abstract}

\maketitle

\section{Introduction}

For the paper \cite{A191} we wanted to apply Theorem 4.11 of Titchmarsh. His theorem is valid for $\sigma\ge\sigma_0>0$, this was sufficient for me, but I needed an explicit constant. 

I searched in the bibliography and found many versions of Theorem 4.11. Most books and papers stated it for $\sigma\ge\sigma_0$, but this hypothesis was not used. Instead many sources used in the proof the bound $|s|\ll t$, even not mentioning this in the statement of the Theorem.

%
%
%
%
The proof by Titchmarsh depends on his Lemma 4.10 which he does not prove. It is also difficult to prove with the conditions given. It is supposed to follow the same lines as the proof of Lemma 4.7, but the monotonicity of the functions needed does not follow from the hypothesis given. In fact, the books by Tenenbaum \cite{Te}*{Cor. 5.1, p.~145} and in Iwaniec and Kowalski \cite{IK}*{eq. (8.28) p.~206} gives a version of Lemma 4.10 by partial summation, giving a different error bound to the one of Titchmarsh and depending on $G$ as in our Lemma \ref{410modified}. 

Iwaniec \cite{Iw}*{Prop. 6.1, p.~24} gives the Theorem, finally, with $\sigma\ge0$ but without explicit constant. The objective of the paper is to get an explicit version of  Proposition 6.1 in Iwaniec \cite{Iw}.

\section{Explicit versions of some lemmas on exponential sums}

\begin{lemma}[Bonnet's form of the Second Mean Value Theorem]
Let $f(x)\in L^1[a,b]$ and $g$ be a real positive and decreasing  function in the range $[a,b]$. Then there is a number $a\le \xi\le b$ such that
\[\int_a^b f(x)g(x)\,dx=g(a)\int_a^\xi f(x)\,dx.\]
\end{lemma}
See \cite{Te}*{Thm.~3, p.~4}.

\begin{lemma}[Explicit version of Lemma 4.3 of Titchmarsh]\label{L:basico2}
Let $f(x)$ and $g(x)$ be continuous real functions defined in a closed interval $[a,b]$. Assume that $f$ has a non-null continuous derivative $f'(x)$ and that 
$g(x)/f'(x)$ is positive and  monotonic. Then 
\[\Bigl|\int_a^b g(x)e^{if(x)}\,dx\Bigr|\le 2\max_{a\le x\le b}\frac{g(x)}{f'(x)}.\]
\end{lemma}
\begin{proof}
For some real $\theta$
\[M:=\Bigl|\int_a^b g(x)e^{i f(x)}\,dx\Bigr|=e^{i\theta}\int_a^b g(x) e^{i g(x)}\,dx=\int_a^b g(x)e^{i (f(x)+\theta)}\,dx\]\[=\Re\int_a^b g(x)e^{i (f(x)+\theta)}\,dx=\int_a^b g(x)\cos(f(x)+\theta)\,dx\]
By the second mean value theorem (assume, for example, that $g(x)/f'(x)$ decreases)
\[M=\int_a^b \frac{g(x)}{f'(x)} f'(x)\cos(f(x)+\theta)\,dx=\frac{g(a)}{f'(a)}\int_a^\xi f'(x)\cos(f(x)+\theta)\,dx\le 2\max_{a\le x\le b}\frac{g(x)}{f'(x)}.\qedhere\]
\end{proof}
\begin{remark}
The constant $2$ in Lemma \ref{L:basico2} is the best possible constant. Take, for example $[a,b]=[0,\pi]$, $g(x)=1$ and $f(x)=x$.  Titchmarsh has a $4$ instead of $2$. Titchmarsh comment: \emph{The values of the constants in these lemmas are usually not of any importance}, which implies that he knew the best constant but did not bother to give it explicitly. But we are interested in giving a good constant.
\end{remark}

\begin{lemma}[Explicit version of Lemma 4.7 of Titchmarsh]\label{L:47modified}
Let $f(x)$ be a real function with a continuous and strictly decreasing derivative $f'(x)$ in $[a,b]$. Let $N\le \lfloor f'(b)\rfloor$ be a non negative integer and define  $\delta=1-(f'(a)-\lfloor f'(a)\rfloor)$
\[\sum_{a<n\le b}e^{2\pi i f(n)}=
\sum_{\nu=N}^{\lfloor f'(a)\rfloor}\int_a^be^{2\pi i(f(x)-\nu x)}\,dx+
\frac{\Theta}{\pi}\Bigl(\pi+3\gamma+3\log(1+f'(a)-N)+\frac{1}{\delta}\Bigr).\]
where $\Theta$ is some complex number with $|\Theta|\le1$ and $\gamma$ is Euler constant.
\end{lemma}
\begin{proof}
We may assume that $N=0$, because our inequality does not change if we take 
$f(x)-Nx$ and $0$ instead of $f(x)$ and $N$. Notice that $\lfloor f'(x)-N\rfloor=
\lfloor f'(x)\rfloor-N\ge0$. 

Hence, we prove our equality assuming that $N=0$. It follows that $N=0\le \lfloor f'(b)\rfloor$, and therefore, $f'(x)\ge f'(b)\ge \lfloor f'(b)\rfloor\ge0$.

Denote by $\psi(x)=\lfloor x\rfloor-x+1/2$, then 
\begin{multline*}\sum_{a<n\le b}e^{2\pi i f(n)}-\int_a^b e^{2\pi i f(x)}\,dx=
\int_a^b e^{2\pi i f(x)}\,d\psi(x)\\=\psi(b)e^{2\pi i f(b)}-\psi(a)e^{2\pi i f(a)}-2\pi i\int_a^b\psi(x)f'(x)e^{2\pi i f(x)}\,dx.\end{multline*}
The function $\psi(x)$ is bounded in absolute value by $1/2$ and has a Fourier series 
\[-\psi(x)=x-\lfloor x\rfloor -\tfrac12=-\frac{1}{\pi}\sum_{\nu=1}^\infty\frac{\sin(2\pi \nu x)}{\nu}\]
with uniformly bounded partial sums. So, we may apply the Dominated Convergence Theorem to interchange sum and integral. It follows that 
\begin{align*}
\sum_{a<n\le b}e^{2\pi i f(n)}-\int_a^b &e^{2\pi i f(x)}\,dx=
2 i \sum_{\nu=1}^\infty \frac{1}{\nu}\int_a^b \sin(2\pi\nu x)f'(x)e^{2\pi i f(x)}\,dx+\Theta,\\
&=\sum_{\nu=1}^\infty\frac{1}{\nu}\int_a^bf'(x)e^{2\pi i(f(x)+\nu x)}\,dx-
\sum_{\nu=1}^\infty\frac{1}{\nu}\int_a^bf'(x)e^{2\pi i(f(x)-\nu x)}\,dx+\Theta
\end{align*}
where $\Theta$ is a complex number of absolute value $\le 1$. 
For $0< \nu\le\lfloor f'(a)\rfloor$   we have 
\begin{multline*}\sum_{\nu=1}^{\lfloor f'(a)\rfloor}\frac{1}{\nu}\int_a^bf'(x)e^{2\pi i(f(x)-\nu x)}\,dx=\sum_{1\le\nu\le\lfloor f'(a)\rfloor}  \frac{1}{\nu}\int_a^b e^{-2\pi i \nu x}d(\tfrac{e^{2\pi i f(x)}}{2\pi i})\\=\sum_{1\le\nu\le\lfloor f'(a)\rfloor} \Bigl.\frac{e^{2\pi i(f(x)-\nu x)}}{2\pi i\nu}\Bigr|_{x=a}^b+\sum_{1\le\nu\le\lfloor f'(a)\rfloor}\int_a^be^{2\pi i (f(x)-\nu x)}\,dx.\end{multline*}
The first term here is bounded in absolute value by 
\[\Bigl|\sum_{1\le\nu\le\lfloor f'(a)\rfloor} \Bigl.\frac{e^{2\pi i(f(x)-\nu x)}}{2\pi i \nu}\Bigr|_{x=a}^b\Bigr|\le \sum_{1\le \nu\le\lfloor f'(a)\rfloor}\frac{1}{\pi\nu}=\frac{\gamma}{\pi}+\frac{\Gamma'(1+\lfloor f'(a)\rfloor)}{\pi\Gamma(1+\lfloor f'(a)\rfloor)}.\]
Hence, we have 
\begin{multline*}\sum_{a<n\le b}e^{2\pi i f(n)}=\sum_{\nu=0}^{\lfloor f'(a)\rfloor}\int_a^be^{2\pi i (f(x)-\nu x)}\,dx+\sum_{\nu=1}^\infty\frac{1}{\nu}\int_a^bf'(x)e^{2\pi i(f(x)+\nu x)}\,dx\\+\sum_{\nu>\lfloor f'(a)\rfloor}\frac{1}{\nu}\int_a^b f'(x)e^{2\pi i (f(x)-\nu x)}\,dx+\frac{\Theta}{\pi}\Bigl(\pi+\gamma+\frac{\Gamma'(1+\lfloor f'(a)\rfloor)}{\Gamma(1+\lfloor f'(a)\rfloor)}\Bigr).
\end{multline*}
We bound the integrals in the second line by means of Lemma \ref{L:basico2}.
Since $\frac{f'(x)}{f'(x)+\nu}$ is decreasing, we have 
\[\Bigl|\int_a^bf'(x)e^{2\pi i(f(x)+\nu x)}\,dx\Bigr|\le \frac{2f'(a)}{2\pi(f'(a)+\nu)}.\]
So,
\[\Bigl|\sum_{\nu=1}^\infty\frac{1}{\nu}\int_a^bf'(x)e^{2\pi i(f(x)+\nu x)}\,dx\Bigr|\le \sum_{\nu=1}^\infty\frac{f'(a)}{\pi\nu(f'(a)+\nu)}=\frac{\gamma}{\pi}+\frac{\Gamma'(f'(a)+1)}{\pi\Gamma(f'(a)+1)}.\]
For $\nu>\lfloor f'(a)\rfloor$, the function $\frac{f'(x)}{2\pi(\nu-f'(x))}\ge0$ and is decreasing, so that by Lemma \ref{L:basico2} we have 
\[\Bigl|\int_a^b f'(x)e^{2\pi i (f(x)-\nu x)}\,dx\Bigr|\le \frac{f'(a)}{\pi(\nu-f'(a))}.\]
Therefore,
\begin{multline*}
\Bigl|\sum_{\nu>\lfloor f'(a)\rfloor}\frac{1}{\nu}\int_a^b f'(x)e^{2\pi i (f(x)-\nu x)}\,dx\Bigr|\\
\le 
\sum_{\nu>\lfloor f'(a)\rfloor}\frac{f'(a)}{\pi\nu(\nu-f'(a))}=\frac{\Gamma'(1+\lfloor f'(a)\rfloor)}{\pi\Gamma(1+\lfloor f'(a)\rfloor)}-\frac{\Gamma'(1+\lfloor f'(a)\rfloor-f'(a))}{\pi\Gamma(1+\lfloor f'(a)\rfloor-f'(a))}.
\end{multline*}
Hence, the three error with the Gamma function are 
\[\frac{\gamma}{\pi}+\frac{\Gamma'(1+\lfloor f'(a)\rfloor)}{\pi\Gamma(1+\lfloor f'(a)\rfloor)}+\frac{\gamma}{\pi}+\frac{\Gamma'(f'(a)+1)}{\pi\Gamma(f'(a)+1)}+\frac{\Gamma'(1+\lfloor f'(a)\rfloor)}{\pi\Gamma(1+\lfloor f'(a)\rfloor)}-\frac{\Gamma'(1+\lfloor f'(a)\rfloor-f'(a))}{\pi\Gamma(1+\lfloor f'(a)\rfloor-f'(a))}\]

We have 
\[\frac{\Gamma'(x)}{\Gamma(x)}\le \log x,\qquad \frac{\Gamma'(x+1)}{\Gamma(x+1)}=
\frac{1}{x}+\frac{\Gamma'(x)}{\Gamma(x)}.\]
It follows that the above terms are 
\[\le \frac{2\gamma+3\log(1+f'(a))}{\pi}+\frac{1}{\pi(1-(f'(a)-\lfloor f'(a)\rfloor))}-\frac{\Gamma'(2+\lfloor f'(a)\rfloor-f'(a))}{\pi\Gamma(2+\lfloor f'(a)\rfloor-f'(a))}\]
Since $1\le 2+\lfloor f'(a)\rfloor-f'(a)\le 2$ and $\Bigl|\frac{\Gamma'(x)}{\Gamma(x)}\Bigr|
\le -\frac{\Gamma'(1)}{\Gamma(1)}=\gamma$ for $1\le x\le 2$,  we get 
\[\le \frac{3\gamma+3\log(1+f'(a))}{\pi}+\frac{1}{\pi\delta}.\]
This ends the proof. 
\end{proof} 

\begin{lemma}[Explicit version of Lemma 4.10 of Titchmarsh]\label{410modified}
Let $f(x)$ be a real function with a continuous and strictly decreasing derivative $f'(x)$ in $[a,b]$. Let $g(x)$ be a continuous function with a continuous derivative. 

Let $N\le \lfloor f'(b)\rfloor$ be a non negative integer and define  $\delta=1-(f'(a)-\lfloor f'(a)\rfloor)$. Then 
\[\sum_{a<n\le b}g(x)e^{2\pi i f(n)}=
\sum_{\nu=N}^{\lfloor f'(a)\rfloor}\int_a^bg(x)e^{2\pi i(f(x)-\nu x)}\,dx+
\frac{\Theta G}{\pi}\Bigl(\pi+3\gamma+3\log(1+f'(a)-N)+\frac{1}{\delta}\Bigr),\]
where $\Theta$ is some complex number with $|\Theta|\le1$ and 
\[G=|g(b)|+\int_a^b|g'(x)|\,dx.\]
\end{lemma}
\begin{proof}
For $a\le y\le b$ we have by Lemma \ref{L:47modified} that 
\[S(y):=\sum_{a<n\le y}e^{2\pi i f(n)}=\sum_{\nu=N}^{\lfloor f'(a)\rfloor}\int_a^ye^{2\pi i(f(x)-\nu x)}\,dx+ R(y)= T(y)+R(y),\]
say,  where 
\[|R(y)|\le R_0:= \frac{1}{\pi}\Bigl(\pi+3\gamma+3\log(1+f'(a)-N)+\frac{1}{\delta}\Bigr).\]
Note the importance here of having  a fixed value of $N$ for all values of $y$, this is possible since we have $N\le \lfloor f'(b)\rfloor\le \lfloor f'(y)\rfloor$. 

Hence, by partial summation (note that $T(a)=0$)
\begin{align*}
\sum_{a<n\le b}g(n)e^{2\pi i f(n)}&=g(b)S(b)-\int_a^b T(y)g'(y)\,dy-\int_a^b R(y)g'(y)\,dy\\
&=S(b)g(b)-T(b)g(b)+\int_a^b g(y)\sum_{\nu=N}^{\lfloor f'(a)\rfloor}e^{2\pi i (f(y)-\nu y)}\,dy-\int_a^b R(y)g'(y)\,dy
\end{align*}
and we have 
\[\Bigl|\int_a^b R(y)g'(y)\,dy\Bigr|\le R_0\int_a^b |g'(y)|\,dy,\qquad 
|S(b)-T(b)|=|R(b)|\le R_0.\]
Hence, we arrive at our result.
\end{proof}

\section{Approximation of  \texorpdfstring{$\zeta(s)$}{zeta(s)}  by a Dirichlet polynomial}

\begin{theorem}[Explicit version of Theorem 4.11 in Titchmarsh]
Let $x>0$ be a given real number and  $s=\sigma+it$ with $\sigma\ge0$, $0<t\le x$, then there is a complex number $\Theta$ with $|\Theta|\le 1$ such that 
\begin{equation}
\zeta(s)=\sum_{n\le x}\frac{1}{n^s}+\frac{x^{1-s}}{s-1}+\frac{29}{14}\frac{\Theta}{x^{\sigma}}.
\end{equation}
\end{theorem}

\begin{proof} Given $x>0$ and  $s=\sigma+it$ with $\sigma>0$ and $0<t\le x$.
Let $M>x$ be a natural number, then we have (for example, see \cite{T}*{eq.~(3.5.3), p.~49})
\[\zeta(s)=\sum_{n\le M}\frac{1}{n^s}+\frac{M^{1-s}}{s-1}-\frac{1}{2M^s}+
s\int_M^\infty \frac{\frac12-\{u\}}{u^{s+1}}\,du.\]
We have 
\[\Bigl|\int_M^\infty \frac{\frac12-\{u\}}{u^{s+1}}\,du\Bigr|\le \frac12\int_M^\infty u^{-\sigma-1}\,du=\frac{1}{2\sigma M^{\sigma}}.\]
Hence,
\[\zeta(s)=\sum_{n\le x}\frac{1}{n^s}+\sum_{x<n\le M}\frac{1}{n^s}+\frac{M^{1-s}}{s-1}+\Theta\Bigl(\frac{1}{2M^\sigma}+\frac{|s|}{2\sigma M^{\sigma}}\Bigr).\]
We now apply Lemma \ref{410modified} to the sum
\[\overline{\sum_{x<n\le M}\frac{1}{n^s}}=\sum_{x<n\le M}\frac{1}{n^\sigma}e^{it\log n}.\]
Hence, we take $g(y)=y^{-\sigma}$, $f(y)=\frac{t}{2\pi}\log y$
with $f'(y)=\frac{t}{2\pi y}$. For $y\in [x,M]$, we have 
\[0<f'(y)=\frac{t}{2\pi y}\le \frac{t}{2\pi x}\le \frac{1}{2\pi}<1.\]
In the lemma, we take $N=0$ and 
\[\delta=1-f'(x)+\lfloor f'(x)\rfloor=1-\frac{t}{2\pi x}\ge 1-\frac{1}{2\pi}.\]
The sum on the right-hand side reduces to the unique value $\nu=0$ so that we get
\[\sum_{x<n\le M}\frac{1}{n^\sigma}e^{it\log n}=\int_x^My^{-\sigma}e^{it\log y}\,dy+\frac{\Theta G}{\pi}\Bigl(\pi+3\gamma+3\log\Bigl(1+\frac{t}{2\pi x}-0\Bigr)+\frac{2\pi}{2\pi-1}\Bigr),\]
where
\[G=M^{-\sigma}+\sigma\int_x^M y^{-\sigma-1}\,dy=M^{-\sigma}+x^{-\sigma}-M^{-\sigma}=x^{-\sigma}.\]
Taking conjugates again
\[\sum_{x<n\le M}\frac{1}{n^s}=\int_x^M y^{-s}\,dy+\frac{\Theta\; x^{-\sigma}}{\pi}\Bigl(\pi+3\gamma+3\log\Bigl(1+\frac{1}{2\pi}\Bigr) +\frac{2\pi}{2\pi-1}\Bigr).\]
Hence, we get 
\begin{multline*}\zeta(s)=\sum_{n\le x}\frac{1}{n^s}+\frac{M^{1-s}}{s-1}+\Theta\Bigl(\frac{1}{2M^\sigma}+\frac{|s|}{2\sigma M^{\sigma}}\Bigr)+\frac{x^{1-s}-M^{1-s}}{s-1}\\+\frac{\Theta\;x^{-\sigma}}{\pi}\Bigl(\pi+3\gamma+3\log\Bigl(1+\frac{1}{2\pi}\Bigr) +\frac{2\pi}{2\pi-1}\Bigr).\end{multline*}
Simplifying and taking limit for $M\to+\infty$ 
\[\zeta(s)=\sum_{n\le x}\frac{1}{n^s}+\frac{x^{1-s}}{s-1}+
\frac{\Theta x^{-\sigma}}{\pi}\Bigl(\pi+3\gamma+3\log(1+\tfrac{1}{2\pi}) +\frac{2\pi}{2\pi-1}\Bigr).\]
The constant is equal to $2.070795<\frac{29}{14}$. Therefore, 
\[\zeta(s)=\sum_{n\le x}\frac{1}{n^s}+\frac{x^{1-s}}{s-1}+
\frac{29}{14}\frac{\Theta}{x^{\sigma}}.\]
Easier to remember
\[\Bigl|\zeta(s)-\sum_{n\le x}\frac{1}{n^s}-\frac{x^{1-s}}{s-1}\Bigr|\le \frac{3}{x^{\sigma}}.\]
For $s=it$ with $t>0$. We take the inequality for  $s'=\sigma+it$ with $\sigma>0$. Taking the limit for $\sigma\to0$ we get the inequality in the case $\sigma=0$. 
\end{proof}

\begin{remark}
As in the book by Titchmarsh, we may give a Theorem for the sum $\sum_{n\le x}n^{-s}$ assuming that $0\le t\le 2\pi x/C$, where $C$ is a given constant $C>1$. 
But we prefer to put $C=2\pi$ to get a concrete version. 
\end{remark}

\end{document}